\documentclass[12equiped]{article}
\usepackage{amsmath,amssymb,amsthm,amscd,mathrsfs}
\usepackage{amsmath,amsfonts,amssymb, amsthm, xypic}
\usepackage[all]{xy}

\input epsf.tex

\title{Remarks on Cohomological Hall algebras and their representations}

\author {Yan Soibelman}

\begin{document}

\maketitle


\newcommand{\CC}{{\mathcal C}}
\newcommand{\LL}{{\mathcal L}}
\newcommand{\MM}{{\mathcal M}}
\newcommand{\NN}{{\mathcal N}}
\newcommand{\OO}{{\mathcal O}}

\renewcommand{\O}{{\mathcal O}}
\newcommand{\Ome}{{\Omega}^{3,0}}
\newcommand{\E}{{\mathcal E}}
\newcommand{\F}{{\mathcal F}}

\newcommand{\g}{{\mathfrak g}}

\renewcommand{\k}{{\bf k}}
\newcommand{\kk}{\overline{\bf k}}

\newcommand{\be}{\begin{equation}}
\newcommand{\ee}{\end{equation}}

\newcommand{\longto}{\longrightarrow}

\newtheorem{thm}{Theorem}[subsection]
\newtheorem{defn}[thm]{Definition}
\newtheorem{lmm}[thm]{Lemma}
\newtheorem{rmk}[thm]{Remark}
\newtheorem{prp}[thm]{Proposition}
\newtheorem{conj}[thm]{Conjecture}
\newtheorem{exa}[thm]{Example}
\newtheorem{cor}[thm]{Corollary}
\newtheorem{que}[thm]{Question}
\newtheorem{ack}{Acknowledgments}

\newcommand{\aA}{{\mathcal A}}
\newcommand{\FF}{{\mathsf{F}}}
\newcommand{\HH}{{\mathcal H}}
\newcommand{\C}{{\bf C}}
\newcommand{\K}{{\bf k}}
\newcommand{\R}{{\bf R}}
\newcommand{\N}{{\bf N}}
\newcommand{\Z}{{\bf Z}}
\newcommand{\Q}{{\bf Q}}
\newcommand{\G}{\mathsf{G}}
\newcommand{\A}{A_{\infty}}
\newcommand{\M}{{\mathsf{M}}}
\newcommand{\epi}{\twoheadrightarrow}
\newcommand{\mono}{\hookrightarrow}
\newcommand\ra{\rightarrow}
\newcommand\uhom{{\underline{Hom}}}
\renewcommand\O{{\cal O}}
\newcommand{\epp}{\varepsilon}

\tableofcontents

\section{Introduction}
The aim of this paper is to  discuss a class of representations of  Cohomological Hall algebras related to the notion of framed stable object of a category. The paper is an extended version of the talk the author gave at the workshop on Donaldson-Thomas invariants at the University Paris-7 in June 2013 and at the conference ``Algebra, Geometry, Physics" dedicated to Maxim Kontsevich (June 2014, IHES). Because of the origin of the paper it contains more speculations than   proofs. 

\subsection{Cohomological Hall algebras and their representations: motivations}

The notion of Cohomological Hall algebra (COHA for short) for quivers with potential was introduced in [KoSo5]\footnote{In fact we considered in the loc.cit. more general case of formally smooth algebras with potential.}. Since a quiver with potential defines a $3$-dimensional Calabi-Yau category ($3CY$ category for short), it was expected that COHA could be defined for ``good"  subcategories of $3$-dimensional Calabi-Yau categories endowed with additional data (most notably,   orientation data introduced in [KoSo1]), or, more generally, for an ind-constructible locally regular $3CY$ category in the sense of [KoSo1]. At the time of writing  this general goal has not been achieved, although there have been some progress in special cases.\footnote{E.g. I heard from Dominic Joyce  about the work in progress of Oren Ben-Bassat in which COHA would be defined for the derived moduli stack of perfect complexes on a compact Calabi-Yau $3$-fold.} 

The definition of COHA given in [KoSo5] is similar to the definition of conventional (constructible) Hall algebra (see e.g. [Sch1]) or its motivic version (see [KoSo1]). Instead of spaces of constructible functions on the stack of objects of an abelian category we considered in [KoSo5] cohomology groups. Certain correspondences define ``structure constants" of the multiplication. There are many versions of COHA  depending on a choice of ``cohomology theory" (see [KoSo5]). It is expected that there is an upgrade of COHA to a dg-algebra in the (triangulated version of the) category of exponential motives (see [KoSo5]). Then all other  versions will appear as ``realizations" of this dg-algebra.

By analogy with conventional Hall algebra of a quiver, which gives a quantization of the ``positive" part of the corresponding Lie algebra, one may want to define ``full" (or ``double") COHA, for which the one defined in [KoSo5] will be just the ``positive part". At this time we do not know the comultiplication which makes COHA into a bialgebra. Hence we cannot follow Drinfeld double construction (which works in the case of  constructible Hall algebras).

On the other hand we can try to define full COHA by means of representation theory, similarly to the classical approach of Nakajima to the infinite Heisenberg algebras (see [Nak2]). In this way one hopes to reconstruct full COHA from its representation theory.

One of the motivations for COHA comes from supersymmetric Quantum Field Theory and String Theory, where spaces of BPS states can be often identified with the cohomology groups of various moduli spaces. From this perspective COHA can be thought of as a mathematical implementation of the idea of BPS algebra (see [HaMo1-2]). Representation theory of BPS algebras has not been developed by physicists, although such a theory should have interesting applications. Furthermore, various dualities in physics can lead to natural mathematical questions about COHA, otherwise unmotivated.
For example, our approach to COHA is based on $3$-dimensional Calabi-Yau categories. The latter appear  in the geometric engineering story on the string-theoretic side. Taking seriously the idea that COHA (or its double) is the BPS algebra, one can ask about the corresponding structures on the gauge-theoretic (i.e.``instanton") side of the geometric engineering. There is some interesting mathematical work related to this question (see e.g. [Nak1], [SchV], [So1], [Sz2]).

In any case the author believes that the representation theory of COHA should be developed further in order to approach some of the above-mentioned problems.

\subsection{COHA and sheaf of vanishing cycles}

Below we will give a very schematic description of COHA assuming that the above-mentioned ``good" category is abelian. 

Let $\k$ be a perfect field.
Suppose $\CC$ is a $\k$-linear triangulated $\A$-category, which is ind-constructible and locally regular in the sense of Section 3 of [KoSo1].\footnote{Calabi-Yau structure which we will discuss later leads to the requirement $char(\k)=0$. For simplicity we will often assume that $\k=\C$.}. It is explained in Section 3.2 of loc.cit that one can associate with $\CC$ the ind-constructible stack $\MM_{\CC}$ of objects of $\CC$. Local regularity implies that $\MM_{\CC}$ is locally presented as an ind-Artin stack over $\k$. Let $\aA\subset \CC$ be an abelian subcategory. Then we have an ind-constructible substack $\MM_{\aA}\subset \MM_{\CC}$ of objects of $\aA$, which is locally ind-Artin.

Hypothetical definition of  COHA  depends on an ind-constructible sheaf $\Phi$ on $\MM_{\CC}$. 
In the case of $3CY$ categories  $\Phi=\phi_W$, which is the sheaf of vanishing cycles of the potential $W$ (we recall the definition of the potential in Section 2.1). For the sheaf of vanishing cycles to be well-defined on $\MM_{\CC}$, the $3CY$-category $\CC$ has to be endowed with an {\it orientation data}. The latter is an ind-constructible super line bundle ${\mathcal L}$ over $\MM_{\CC}$, such that for the fiber over a point of $\MM_{\CC}$ corresponding to an object  $E\in Ob(\CC)$ one has ${\mathcal L}_E^{\otimes 2}=sdet(Ext^{\bullet}(E,E))$. Furthermore, it is required that  
${\mathcal L}_E$ behaves naturally on exact triangles (see [KoSo1], Section 5 for the details). It follows from local regularity that ${\mathcal L}$ is (locally) a line bundle over an ind-Artin stack.

Let $i:\MM_{\aA}\subset \MM_{\CC}$ be the natural embedding. Then the pull-back $i^{\ast}(\Phi)$ is an ind-constructible sheaf on $\MM_{\aA}$.
Let ${\mathcal Z}\subset \MM_{\aA}\times \MM_{\aA}$ be the ``Hecke correspondence", which is the stack consisting of pairs of objects $(E,F)$ such that $E\subset F$. There are projections $p_n: {\mathcal Z}\to \MM_{\aA}, n=1,2,3$ such that $p_1(E,F)=E, p_2(E,F)=F, p_3(E,F)=E/F$. One of the properties  that we require from ``good" abelian category $\aA$ is  that the projection $p_1$ is a proper morphism of ind-Artin stacks.

As a vector space COHA of $\MM_{\aA}$ is defined as  ${\mathcal H}:={\mathcal H}_{\aA}={\bf H}^{\bullet}(\MM_{\aA}, i^{\ast}(\Phi))$. For that one chooses an appropriate cohomology theory ${\bf H}^{\bullet}$ of Artin stacks with coefficients in constructible sheaves.
The product $m: {\mathcal H}_{\aA}\otimes {\mathcal H}_{\aA} \to {\mathcal H}_{\aA}$ is defined by the formula
$p_{1\ast}\circ (p_2^{\ast}\otimes p_3^{\ast})(i^{\ast}(\Phi)\otimes i^{\ast}(\Phi))$. 

For the sheaf of vanishing cycles $\phi_W$ associativity of the product depends on  the Thom-Sebastiani theorem for the chosen cohomology theory.\footnote{ As explained in Section 7 of [KoSo5], it is more convenient to work with compactly supported cohomology and then apply the duality functor.}

\begin{rmk}
As a part of the data we fix an abstract version of Chern character (called the class map in [KoSo1]). It is a homomorphism of abelian groups $cl: K_0(\CC)\to \Gamma$, where $\Gamma\simeq \Z^n$ is a free abelian group, such that connected components of  $\MM_{\CC}$ are parametrized by $\Gamma$, while classes $cl(E), E\in Ob(\aA)$ form an additive submonoid $\Gamma_+\subset \Gamma$. Then COHA of $\aA$ will be $\Gamma_+$-graded algebra.
\end{rmk}
In the case of smooth algebras with potential considered in [KoSo5] the stack $\MM_{\aA}$ is a countable union of smooth quotient stacks, and the foundational questions are resolved positively. For some ideas about the general case one can look at [DyKap]. Natural class of examples for which COHA should be defined comes from the theory of graded symplectic manifolds (see [PaToVaVe] about foundational material and [BuJoMe] about some progress in the  definition of $\phi_W$).

\subsection{Stable framed objects and modules over COHA}

It is natural to ask whether one can realize representations of COHA in the cohomology groups of some natural schemes (or stacks), which might also depend on a choice of stability condition on $\aA$.
Let us explain how it can be achieved.\footnote{We warn the reader that even in the case of quivers our moduli spaces are {\bf not}  Nakajima quiver varieties.}

First we define the moduli space $\MM_{\gamma}^{fr,st}, \gamma \in \Gamma_+$ of `` stable framed objects of class $\gamma$" (in applications those can be framed stable sheaves, framed stable representations of quivers, framed special Lagrangian submanifolds, etc.).
This notion depends on a choice of stability condition on $\aA$. It is expected (see [KoSo7]) that the  sheaf $\phi_W$ ``descends" to each moduli space $\MM_{\gamma}^{fr,st}$.

For a pair of classes $\gamma_1,\gamma_2\in \Gamma_+$  let us consider the Hecke correspondence ${\mathcal Z}_{\gamma_1,\gamma_2}$ of pairs $(E_{\gamma_1+\gamma_2}, E_{\gamma_2})$ (the subscripts denote the Chern classes) of framed stable objects such that $E_{\gamma_2}$ is a quotient of
$E_{\gamma_1+\gamma_2}$. Let us denote the  cohomology theory we used in the definition of COHA by ${\bf H}^{\bullet}$. It descends to each $\MM_{\gamma}^{fr,st}$. Furthermore, similarly to the definition of COHA we have three projections of  ${\mathcal Z}_{\gamma_1,\gamma_2}$: 

a) to  $\MM_{\gamma_2}^{fr,st}$;

b) to the moduli space $\MM_{\gamma_1}$ of all (not framed) objects with fixed $\gamma_1$;

c) to  $\MM_{\gamma_1+\gamma_2}^{fr,st}$.

Using the pull-back and pushforward construction as in the previous subsection, we obtain a structure of ${\mathcal H}_{\aA}=\oplus_{\gamma}{\bf H}^{\bullet}(\MM_{\gamma},\phi_W)$-module over COHA of $\aA$ on the space $\oplus_{\gamma}{\bf H}^{\bullet}(\MM_{\gamma}^{fr,st},\phi_W)$.

We are going to discuss stable framed objects in Section 3. We will show there that the moduli stacks of stable framed objects are in fact schemes. Hence graded components of our representations of COHA are {\it finite-dimensional} vector spaces. As a version of the above considerations we can drop the stability assumption and consider stacks of framed (but not necessarily stable framed) objects. Then we still obtain  representations of COHA.  But this time the graded components  of the representation spaces will be in general infinite-dimensional.

\subsection{Modules over COHA motivated by physics}

 As we have already mentioned, COHA can be thought of as a mathematical incarnation of the  algebra of (closed) BPS states envisioned in [HaMo1,2]. According to [KoSo5] (refined) BPS-invariants  can be computed as virtual Poincar\'e polynomials of graded components of COHA. One can ask about the meaning of the algebra structure on COHA.
 
 Motivated by the ideas of S. Gukov (see e.g. [GuSto]) we would like to think about representations of COHA described in the previous subsection as of  representations of the algebra of closed BPS states on the vector space of open BPS states. We are going to speculate about applications of this point of view in the last section of the paper. We plan to discuss the relationship between COHA and BPS-algebras more systematically in separate projects jointly with E. Diaconescu, S.Gukov, N. Saulina.

Here we just mention three interesting classes of representations of  COHA which have geometric origin and should have interesting applications to gauge theory and knot invariants:

a) Representation of COHA of the resolved conifold $X={tot(\mathcal O}_{{\bf P}^1}(-1)\oplus {\mathcal O}_{{\bf P}^1}(-1))$ realized in the cohomology of moduli spaces of $C$-framed stable sheaves in the sense of [DiHuSo]. Those modules should have applications in the theory of algebraic knots and Hilbert schemes of curves (see [ORS], [DiHuSo], [Mau1]).

b) Modules over COHAs of the Fukaya categories of non-compact Calabi-Yau $3$-folds associated with spectral curves of Hitchin integrable systems. Those should serve as BPS algebras of some gauge theories from the class ${\mathcal S}$ (see e.g. [GaMoNe-2], [Ga1] about the latter).

c) This class of examples was already mentioned above. Based on the ideas of geometric engineering we hope for a class of representations of  COHA related to the moduli spaces of (framed, possibly ramified) instantons on ${\bf P}^2$.  One can hope to understand the relationship between the algebras of Hecke operators proposed in [Nak1] and those proposed in [So1].  Currently  we can see  COHA on the ``Calabi-Yau side" of geometric engineering and the (seemingly unrelated) conventional (``motivic") Hall algebra on the ``instanton side" (cf. also [SchV], [Sz2]). 

The relationship between various classes of gauge theories might give non-trivial results about corresponding COHAs and their representations
(including the relations between  a), b), c)).

\subsection{Contents of the paper}

Section 2 is a reminder on COHA in the framework of quivers with potential. Section 3 is devoted to    stable framed objects in triangulated and abelian categories. In Section 4 we discuss representations of COHA realized in the cohomology of the moduli spaces of stable framed  representations.  We also discuss an approach to the definition of  ``full COHA" based on the representation theory in the cohomology of moduli of stable framed objects. In Section 5 we speculate about representations of COHA motivated by knot theory and physics.

 {\it Acknowledgments.} I thank to Alexander Braverman, Emanuel Diaconescu, Sergei Gukov, Nigel Hitchin, Mikhail Kapranov, Gabriel Kerr, Melissa Liu, Hiraku Nakajima, Andy Neitzke, Alexei Oblomkov, Vivek Schende, Olivier Schiffmann, Xinli Xiao for stimulating discussions and correspondences. I am especially grateful to Maxim Kontsevich for numerous discussions on the subject of this paper. I also thank IHES for excellent research conditions. This work is partially supported by an NSF grant.

\section{Cohomological Hall algebra: reminder}

This section is a reminder of some basic facts about the notion of Cohomological Hall algebra. 
Most of the material is borrowed from [KoSo5], and we refer the reader to loc.cit. for more details and proofs.

\subsection{COHA and $3CY$ categories}

 Suppose we are given an ind-constructible locally regular $3CY$ category $\CC$ over the field $\k, char(\k)=0$ (see [KoSo1]). As explained in Section 3.2 of loc.cit., one can associate with such a category  the stack  of objects, which is a countable disjoint union of schemes over $\k$ of finite type acted by affine algebraic groups. For simplicity of the exposition we take the ground field $\k=\C$.

Some examples of such categories are listed in the Introduction of  [KoSo1]. They include various categories of $D$-branes popular in string theory (e.g. the Fukaya category of a compact or local Calabi-Yau $3$-fold, the category of perfect sheaves on such a $3$-fold, the category of finite-dimensional representations of a quiver with potential, etc.).

In order to define COHA one  has  to choose   {\it orientation data} (see [KoSo1], Section 5) on $\CC$ as well as a ``good" $t$-structure with the ind-Artin heart. Let us denote it by ${\mathcal A}$. The existence of mutation-invariant orientation data is known for a class of $3CY$ categories associated with a quiver without potential (see [Dav1]). There are partial existence results for the derived category of coherent sheaves on a compact Calabi-Yau $3$-fold (see e.g. [Hu]). Probably Dominic Joyce with collaborators will construct soon an orientation data on the stack of objects of the category of perfect complexes on a Calabi-Yau $3$-fold. But the general case is still open. In present  paper we will assume the existence of the orientation data as a part of the ``foundational" package. Also, we do not discuss in detail the meaning of the notion of ``good" $t$-structure. As we mentioned in the Introduction, the latter includes properness of the morphisms which appear in the definition of the product on COHA.

We assume as part of the data the ``class map" $cl:K_0(\CC)\to \Gamma$ (see [KoSo1]), where $\Gamma\simeq \Z^n$  is a free abelian group endowed with integer skew-symmetric form $\langle\bullet,\bullet\rangle$ (Poisson lattice). We also assume that the class map respects  the Euler form $\chi(E,F)=\sum_i(-1)^idim\,Ext^i(E,F)$ on $K_0(\CC)$ and the form $\langle\bullet,\bullet\rangle$ on $\Gamma$. The lattice $\Gamma$ plays a role of topological $K$-theory of the category $\CC$. Finally, we assume that we have fixed an additive submonoid $\Gamma_+\subset \Gamma$  generated by $cl(E), E\in Ob({\mathcal A})$. 

When the above choices are made, one can define COHA of ${\mathcal A}$ as an associative algebra graded by $\Gamma_+$. Graded components are given by the cohomology of the moduli stacks of objects with the given  class $\gamma\in \Gamma$ with the coefficients in the sheaf of vanishing cycles of the potential of $\CC$ restricted to ${\mathcal A}$. 

For completeness we recall here the notion of  potential of a $3CY$ category.
Using the $\A$-structure on  $\CC$ as well as the Calabi-Yau pairing $(\bullet,\bullet)$ (see [KoSo1]) one defines the   potential of an object $E$ as a formal series:

$$W_E(a)=\sum_{n\ge 1}{\left(m_n(a,...,a),a\right)\over{n+1}},$$
where $m_n$ are higher composition maps, and the element $a$ belongs to $Hom^1(E,E)$ which is the subspace in the graded space $Hom(E,E)$ consisting of elements of degree $1$. By our assumptions the potential $W_E$ is a locally regular function with respect to $E$. Hence we have a partially formal function $W$ defined by the family of series $W_E$.

\begin{rmk}
If $\CC$ is ``minimal on the diagonal" (see [KoSo1]), we can replace $Hom(E,E)$ by its cohomology with respect to the differential $m_1$. In this case we may assume that $a\in Ext^1(E,E)$, which can be thought of as the ``tangent space to the moduli stack of formal deformations of $E$". Hence one can think of the potential as a  function on the moduli stack of objects which is locally regular along the stack of objects (this follows from the ``locally regular" assumption) and formal in the transversal direction.
\end{rmk}

Then COHA is a $\Gamma$-graded vector space
$$\mathcal{H}:=\oplus_{\gamma\in \Gamma} \mathcal{H}_\gamma\,\,,$$
where $\mathcal{H}_\gamma=H_{\G_{\gamma}}^{\bullet}(S_{\gamma}, W_{\gamma})$, and $S_{\gamma}$ is the stack of objects $E$ such that $cl(E)=\gamma$.
Recall that we  use  an appropriate stack version of the cohomology theory $H^{\bullet}(X,f)$ of a scheme $X$ endowed with a regular function $f$. There are several choices for such  theory. They are discussed in [KoSo5], where the above approach made rigorous in the case of $3CY$ categories arising from quivers (more generally, formally smooth algebras) with potential. A version of the cohomology theory which is suitable in the framework of categories is called ``critical cohomology" in loc. cit. It is defined by means of the compactly supported cohomology of $X$ with coefficients in the sheaf of vanishing cycles of $f$. Sometimes (e.g. for quivers with potential) the function $f:=W$ is regular. In such a case one can use de Rham cohomology defined via the twisted de Rham differential $d+dW\wedge(\bullet)$ or Betti cohomology which is generated by ``integration cycles" for the exponential differential forms of the type $exp(W)\nu$. More generally, one can define ``motivic" version of COHA. In that case COHA $\mathcal{H}$ is an object of the tensor category of exponential mixed Hodge structures, and the concrete choice of the cohomology theory corresponds to a tensor functor to graded vector spaces (``realization").
It is explained in [KoSo5] that in all realizations $\mathcal{H}$ carries an associative algebra structure with ``structure constants" defined by means of the cohomology of certain ``Hecke correspondences" with coefficients in the sheaves of vanishing cycles of the potential $W=(W_{\gamma})_{\gamma\in \Gamma}$.

Let us illustrate the above considerations in the case of a quiver $Q$ with potential $W$, which is the main example in [KoSo5]. We set $\k=\C$. If $I$ is the set of vertices of $Q$ then $\Gamma=\Z^I, \Gamma_+=\Z_{\ge 0}^I$. For any $\gamma=(\gamma^i)_{i\in I}\in \Gamma_+$ we consider $\gamma$-dimensional representations of $Q$ in coordinate vector spaces $(\C^{\gamma^i})_{i\in I}$. It is an affine scheme $\M_{\gamma}$ naturally acted by the affine algebraic group $\G_{\gamma}=\prod_{i\in I}GL(\gamma^i,\C)$. Then the corresponding stack of objects is a countable union (over all dimension vectors $\gamma\in \Gamma_+$)  of algebraic varieties $Crit(W_{\gamma})$ of the critical points of the functions $W_{\gamma}=Tr(W):\M_{\gamma}\to \C$. Then COHA is the direct sum $\oplus_{\gamma\in \Gamma_+}H^{\bullet}_{\G_{\gamma}}(\M_{\gamma},W_{\gamma})$ with the product defined in the loc.cit. In the next three subsections we are going to recall more explicit descriptions of the product in some examples.

\subsection{COHA for quivers without potential}

COHA is non-trivial even if $W=0$. In the latter case
$$\mathcal{H}_\gamma:=H^\bullet_{\G_\gamma} (\M_\gamma).$$
 Since $\M_{\gamma}$ is equivariantly contractible, and $\G_{\gamma}$ is homotopy equivalent to its maximal torus, one can use the toric localization and obtain an explicit formula for the product which expresses COHA as a shuffle algebra. In the formula below we identify equivariant cohomology of a point with respect to the trivial action of the torus $(\C^{\ast})^n$ with the space of symmetric polynomials in $n$ variables.

\begin{thm}
The product $f_1\cdot f_2$ of elements $f_i\in \mathcal{H}_{\gamma_i},\,i=1,2$  is given by the symmetric function
$g((x_{i,\alpha})_{i\in I, \alpha\in \{1,\dots,\gamma^i\}})$, where $\gamma:=\gamma_1+\gamma_2$,  obtained from
the following function in variables $({x}'_{i,\alpha})_{i\in I, \alpha\in \{1,\dots,\gamma^i_1\}}$ and $({{x}}''_{i,\alpha})_{i\in I, \alpha\in \{1,\dots,\gamma^i_2\}}$:

$$f_1(({x}'_{i,\alpha}))\,
f_2(({x}''_{i,\alpha})) \,\,
\frac{ \prod_{i,j\in I}\prod_{\alpha_1=1}^{\gamma^i_1}\prod_{\alpha_2=1}^{\gamma^j_2} ({{x}}''_{j,\alpha_2}
-{x}'_{i,\alpha_1}  )^{a_{ij}
}}{\prod_{i\in I}
\prod_{\alpha_1=1}^{\gamma^i_1}\prod_{\alpha_2=1}^{\gamma^i_2}({{x}}''_{i,\alpha_2}-{x}'_{i,\alpha_1})
}\,\,,$$
by taking the sum over all shuffles for any given $i\in I$ of the variables
${x}'_{i,\alpha},{{x}}''_{i,\alpha}$ (the sum is over $\prod_{i\in I}\binom{\gamma^i}{\gamma_1^i}$ shuffles).

\end{thm}

Here $a_{ij}$ is the number of arrows in $Q$ from the vertex $i$ to vertex $j$.

For example, let $Q=Q_d$ be  a quiver with just one vertex and $d\ge  0$ loops. Then the product formula  specializes
to
$$(f_1\cdot f_2)(x_1,\dots,x_{n+m}):=$$
$$\sum_{i_1,...,j_m} f_1(x_{i_1},\dots,x_{i_n})\,f_2(x_{j_1},\dots,x_{j_m})\,
\left(\prod_{k=1}^n\prod_{l=1}^m(x_{j_l}-x_{i_k})\right)^{d-1}
$$
for symmetric polynomials, where $f_1$ has $n$  variables, and  $f_2$ has $m$ variables.
The sum is taken over all $\{i_1<\dots<i_n, j_1<\dots<j_m,
\{i_1,\dots,i_n,j_1,\dots,j_m\}
=\{1,\dots,n+m\}$.
The product $f_1\cdot f_2$ is a symmetric polynomial
in $n+m$ variables. One can show that for even $d$ the algebra is isomorphic to the infinite Grassmann algebra, while for odd $d$ one gets an infinite symmetric algebra.

We introduce a double grading on  algebra ${\cal H}$, by declaring that a homogeneous symmetric polynomial of degree $k$ in $n$ variables
has bigrading $(n,2k+(1-d)n^2)$. Equivalently, one can shift the cohomological grading in $H^\bullet(\mathrm{BGL}(n,\C))$ by
$[(d-1)n^2]$. In  general, even for quivers without potential  each component $\mathcal{H}_\gamma$ has also the grading by cohomological degree. Total $\Gamma\times \Z$-grading can be further refined, since $\mathcal{H}_\gamma$ carries the weight filtration (as an object of the category of exponential mixed Hodge structures, see [KoSo5]). Hence typically COHA has $\Gamma\times \Z\times \Z$-grading (which is not compatible with the product). More precisely, it is shown in [KoSo5] that for $W=0$ COHA is graded by the Heisenberg group.

Finally, we remark that in the case of Dynkin quivers there are other interesting explicit formulas for the product in COHA (see [Rim]).

\subsection{COHA for quiver $A_2$}

The quiver $A_2$ has two vertices $\{1,2\}$ and one arrow $1\leftarrow 2$. The Cohomological Hall algebra $\mathcal{H}$ of this quiver contains two subalgebras $\mathcal{H}_L,\,\mathcal{H}_R$ corresponding to representations supported at the vertices $1$ and $2$ respectively.
Clearly each subalgebra $\mathcal{H}_L,\,\mathcal{H}_R $ is isomorphic to the  Cohomological Hall algebra for  the quiver $A_1=Q_0$. Hence it is
an infinite Grassmann algebra. Let us denote the generators by $\xi_i,\, i=0,1,\dots$ for the vertex $1$ and by $\eta_i,\, i=0,1,\dots$ for the vertex $2$. Each generator $\xi_i$ or $\eta_i$ corresponds to an additive generator of the group $H^{2i}(BGL(1,\mathbb{C}))\simeq \mathbb{Z}\cdot x^i$. Then one can check that  $\xi_i,\eta_j, \,i,j\geqslant 0$ satisfy the relations

$$\xi_i\xi_j+\xi_j\xi_i=\eta_i\eta_j+\eta_j\eta_i=0\,,\,\,\,\,\eta_i\,\xi_j=\xi_{j+1}\eta_i-\xi_j\eta_{i+1}\, .$$

Let us introduce the elements $\nu_i^1=\xi_0\eta_i\,,\, i\geqslant 0$ and $\nu_i^2=\xi_i\eta_0\,, \,i\geqslant 0$. It is easy to see that $\nu_i^1\nu_j^1+\nu_j^1\nu_i^1=0$, and similarly the generators $\nu_i^2$ anticommute. Thus we have two infinite Grassmann subalgebras in $\mathcal{H}$ corresponding to these two choices: $\mathcal{H}^{(1)}\simeq \bigwedge(\nu_i^1)_{i\geqslant 0}$ and

$\mathcal{H}^{(2)}\simeq \bigwedge(\nu_i^2)_{i\geqslant 0}$. One can directly check the following result.

\begin{prp} The multiplication (from the left to the right) induces   isomorphisms of graded vector spaces:

$$\mathcal{H}_L\otimes \mathcal{H}_R\stackrel{\sim}{\longrightarrow}\mathcal{H}, \,\,\,\,\,\,
\mathcal{H}_R\otimes \mathcal{H}^{(i)}\otimes \mathcal{H}_L\stackrel{\sim}{\longrightarrow} \mathcal{H}\,,\,\,i=1,2\,. $$

\end{prp}

\subsection{COHA for Jordan quiver  with polynomial potential}

Let us consider  the quiver $Q_1$ which has one vertex and one loop $l$ (Jordan quiver), and choose  as the potential $W=\sum_{i=0}^{N} c_i l^i,\,\,\,c_N\ne 0$ an arbitrary polynomial of degree $N\in \Z_{\geqslant 0}$ in one variable.

In the case $N=0$, the question about COHA reduces to the quiver $Q_1$ without potential. This case was considered before. The algebra $\mathcal{H}$ is the symmetric algebra of infinitely many variables. 

In the case $N=1$ COHA is one-dimensional.

In the case $N=2$ we may assume without loss of generality that $W=-l^2$. Then COHA $\mathcal{H}=\mathcal{H}^{(Q_1,W)}$ is the exterior algebra with infinitely many generators (infinite Grassmann algebra). This can be shown directly.

In the case when the degree $N\geqslant 3$, one can show that the bigraded algebra $\mathcal{H}$
is isomorphic to the $(N-1)$-st tensor power of the  infinite Grassmann algebra of the case $N=2$.

Basically the above examples are the only cases in which we know COHA explicitly. On the other hand, generating functions for the dimensions of its graded components (we call them {\it motivic DT-series} in [KoSo1,5]) are known in many cases.

\subsection{Stability conditions and motivic DT-invariants}

Definition of COHA depends on the abelian category $\aA$ but does not depend on the central charge, which is a homomorphism of groups $Z:\Gamma\to \C$. This raises the question about the role of Bridgeland stability condition in the structure of COHA.

Having a central charge $Z:\Gamma\to \C$ we can define a full subcategory ${\mathcal A}$ of our category $\CC$ generated by zero object and semistable objects with the central charge sitting in a given strict sector $V\subset \R^2$ (the sector  has the vertex in the origin). For example, we can take $V=l$ to be a ray. Taking $V$ to be the upper-half plane we arrive to a category which is the heart of a $t$-structure of our $3CY$ category $\CC$. In these two cases the categories generated by semistables are abelian.

 As explained in [KoSo5], for a fixed strict sector $V$, one can define a $\Gamma$-graded {\it vector space}
$$\mathcal{H}(V):=\oplus_{\gamma\in \Gamma} \mathcal{H}_\gamma(V)\,\,.$$
But this space cannot be endowed  with a structure of an associative algebra, except of the case when $V=l$ or $V$ being an upper-half plane. The problem is with properness of morphisms of the corresponding stacks.

It was observed in [KoSo5], Section 5.2 that the algebras $\mathcal{H}_l:=\mathcal{H}(l)$ resemble universal enveloping algebras of some Lie algebras $\g_l$ which are analogous to the ``positive root" Lie algebras $\g_{\alpha},\alpha>0$ of Kac-Moody algebras. 
Then similarly to the isomorphism $U(n_+)\simeq \otimes_{\alpha>0}U(\g_{\alpha})$ (which depends on a chosen order on the set of positive roots) one should expect an isomorphism $\mathcal{H}(V)\simeq\otimes_{l\subset V}\mathcal{H}_l$ where the tensor product is taken in the clockwise order over all rays in the sector $V$. This was demonstrated in [Rim] in the case of Dynkin quivers without potential. In particular, taking $V$ to be the upper-half plane we obtain a factorization of the  COHA $\mathcal{H}$ into the tensor product of COHAs  for individual rays. COHA for each ray $l$ is typically commutative. It can be computed from the knowledge of space of semistable objects in the fixed $t$-structure whose central charges belong to $l$. For a generic central charge we have two possibilities: either $l$ does not contain $Z(\gamma)$ for $\gamma\in \Gamma$, or $l$ contains only multiples $nZ(\gamma_0), n>0$ for some primitive vector $\gamma_0$ (an furthermore, only vectors $n\gamma_0, n\in \Z_{>0}$ are mapped by $Z$ to $l$). In this case $\mathcal{H}_l$ is indeed commutative and can be computed explicitly in many cases.

The notion of motivic DT-series (i.e. virtual Poincar\'e series of $\mathcal{H}$) does not depend on the central charge.
On the other hand, motivic {\it DT-invariants} $\Omega^{mot}(\gamma)$ (they correspond in physics to refined BPS invariants) can be defined only after a choice of  stability condition (i.e. the central charge in case of quivers). Definition of DT-invariants is based on the theory of factorization systems developed in [KoSo5]. It follows from loc. cit. that the motivic DT-series factorizes as a product of the powers of shifted quantum dilogarithms. Those powers are motivic DT-invariants.

As a side remark we mention that factorization systems appear in different disguises when mathematicians try to make sense of the operator product expansion in physics (we can mention e.g. the work of Beilinson and Drinfeld on chiral algebras or the work of Costello and others on OPE in QFT). From this point of view it is not quite clear why factorization systems appear in our story.

\subsection{Generators of COHA}

For different $t$-structures the corresponding COHAs are not necessarily isomorphic. For example, if we start with a pair $(Q,W)$ consisting of a quiver $Q$ with potential $W$ and make a mutation at a vertex $i_0\in I$, then COHA for the mutated pair $(Q^{\prime}, W^{\prime})$ is different from the one for $(Q,W)$. On the other hand we can compute motivic DT-series for the mutated quiver with potential. As was explained in [KoSo1] and [KoSo5], if we make a mutation at the vertex $i_0\in I$ then the motivic DT-series for $(Q,W)$ and $(Q^{\prime}, W^{\prime})$ are related by the conjugation by the motivic DT-series corresponding to the ray $l_0=\R_{>0}\cdot Z(\gamma_{i_0})$ (which is essentially the quantum dilogarithm).

\begin{que} How to define COHA for a {\it triangulated} $3CY$ category $\CC$? 
\end{que}
We do not know the answer to this question, but we can see some structures which should be incorporated in the definition.

For example, let us consider all COHAs corresponding to all possible mutations. Let $M$ be the orbit of the pair $(Q,W)$ under the action of the group of mutations. Then to any $m\in M$ we can assign COHA $\mathcal{H}_m$.
More generally, we can consider rotations $Z\mapsto Ze^{i\theta}$ of the central charge and get the corresponding COHA $\mathcal{H}_{e^{i\theta}}$. This defines a structure of cosheaf of algebras over $S^1$. Each stalk is the COHA 
for the corresponding $t$-structure.

Next question is about the space of generators of COHA.
Recall the following conjecture from [KoSo5] which was proved by Efimov (see [Ef]). It is formulated for symmetric quivers. Such quivers arise naturally in relation to $2$-dimensional Calabi-Yau categories and Kac-Moody algebras.

\begin{thm} Let ${\mathcal H}$ be the COHA (considered as an algebra over $\Q$) for the abelian category of finite-dimensional representations of a symmetric quiver $Q$. Then ${\mathcal H}$ is a free supercommutative algebra generated by a graded vector space $V$ over $\Q$ of the form
$V=V^{\prime}
\otimes \Q[x]$, where $x$ is an even variable of
bidegree $(0,2)\in \Z_{\ge 0}^I\times \Z$, and for any given $\gamma$ the space $V^{\prime}_{\gamma,k}\ne 0$ is non-zero (and finite-dimensional) only for {\bf finitely many} $k\in \Z$.
\end{thm}

In general we  expect (see [KoSo5] for the precise question) that $\mathcal{H}$ is isomorphic to the universal enveloping algebra of a graded Lie algebra $V:=V^{\prime}
\otimes \C[x]$ which satisfies the conditions of the Theorem 2.6.2. Mutations act on $V$, hence we obtain a collection of vector spaces $V_m$ (one for each $t$-structure $m$). From the point of view of chamber structure of the space of stability conditions, we can say that with every chamber we associate its own COHA. Change of the chamber corresponds to the wall-crossing, which at the level of COHA is a conjugation (with a shift of grading). More structural results generalizing on COHA, including generalizations of the above Theorem 2.6.2 and their applications (e.g. to Kac conjecture) can be found in recent papers by B. Davison (see [Dav2,3]).

\section{Framed and stable framed objects}

In this section we present a definition of  stable framed objects following [KoSo7] as well as a related construction of  modules over COHA of the same authors (unpublished).

\subsection{Stable framed objects in triangulated  categories}


We recall the definition of stable framed object from [KoSo7] in the case of triangulated categories. Then we discuss some versions in the case of abelian categories.

Let $\CC$ be a triangulated $\A$-category over the ground field $\k$, which we assume to be an algebraically closed of characteristic zero. We  fix a stability condition $\tau\in Stab(\CC)$. Let  $\Phi:\CC\to D^b(Vect_{\k})$ be an exact functor to the triangulated category of bounded complexes of $\k$-vector spaces.

 For a fixed ray $l$  in the upper-half plane with the vertex at the origin, we denote by $\CC_{l}:=\CC^{ss}_{l}$ the abelian category of $\tau$-semistable objects having  the central charge in $l$. We will impose the following assumption: {\it $\Phi$ maps  $\CC_{l}$ to the complexes concentrated in non-negative degrees}.

\begin{defn} Framed object (or $\Phi$-framed object, if we want to stress dependence on the framing functor) is a pair $(E,f)$ where $E\in Ob(\CC_{l})$ and $f\in H^0(\Phi(E))$.

\end{defn}

Let $(E_1,f_1)$ and $(E_2,f_2)$ be two framed objects. We define a morphism $\phi:(E_1,f_1)\to (E_2,f_2)$ as a morphism $E_1\to E_2$ such that the induced map $H^0(\Phi(E_1))\to H^0(\Phi(E_2))$ maps $f_1$ to $f_2$. Framed objects naturally form a category, and hence there is a notion of isomorphic framed objects.

\begin{defn} We call the framed object $(E,f)$ stable is there is no  exact triangle $E^{\prime}\to E\to E^{\prime\prime}$ in $\CC$ with $E^{\prime}$ non-isomorphic to $E$ such that both $E^{\prime},E^{\prime\prime}\in Ob(\CC_{l})$ and such that there is
$f^{\prime}\in H^0(\Phi(E^{\prime}))$ which is mapped to $f\in H^0(\Phi(E))$.

\end{defn}

Then one deduces the following result 
(see [KoSo7]), proof of which we reproduce here for completeness.                                              

\begin{prp} If $(E,f)$ is a stable framed object then $Aut(E,f)=\{1\}$.

\end{prp}

{\it Proof.} Let $h\in Aut(E)$ satisfies the property that its image $\Phi(h)$ preserves $f$. We may assume that $h\in Hom^0(E,E)$.  We would like to prove that $h=id$. Assume the contrary. Let $h_1:=h-id\ne 0$. Then $\Phi(h_1)(f)=0$. Since the category $\CC_l$ is abelian, the morphism $h_1\ne 0$ gives rise to a  short exact sequence in $\CC_l$:

$$0\to Ker(h_1)\to E\to Im(h_1)\to 0,$$
where $Im(h_1)\ne 0$. Hence there exists an exact triangle $E^{\prime}\to E\to E^{\prime\prime}$  in $\CC$ with $E^{\prime}=Ker(h_1)$ non-isomorphic to $E$ and $E^{\prime\prime}=Im(h_1)$. Let us consider a short exact sequence in $\CC_l$ given by
$$0\to Ker(h_1)\to E\to E\to Coker(h_1)\to 0,$$
where the morphism $E\to E$ is $h_1$.
Since the functor $\Phi$ is exact we get short exact sequence of vector spaces

$$H^0(Ker(\Phi(h_1)))\to H^0(\Phi(E))\to H^0(\Phi(E))\to H^0(Ker(\Phi(h_1)))\to H^1(Ker(\Phi(h_1)))\to....$$

Let us remark that by the assumption that $\Phi$ maps $\CC_l$ to complexes with non-negative cohomology, we conclude that if $E^{\prime}\to E\to E^{\prime\prime}$  is an exact triangle then in the induced exact sequence
$$H^{-1}(\Phi(E^{\prime\prime}))\to H^0(\Phi(E^{\prime}))\to H^0(\Phi(E))\to...$$
the first terms is trivial.
Hence the functor $H^0\Phi$ maps monomorphisms in $\CC_l$ to monomorphisms in the category $Vect_{\k}$ of $\k$-vector spaces.

Let us decompose $h_1$ into a composition of  the morphism $\psi: E\to Im(h_1)$ and the natural embedding $j:Im(h_1)\to E$. Applying $\Phi$, and using  $\Phi(h_1)(f)=0$ and the above remark we conclude that 
$\Phi(\psi)(f)=0$.

Finally, applying $\Phi$ to the short exact sequence
$$0\to Ker(h_1)\to E\to Im(h_1)\to 0,$$
we obtain a short exact sequence in $Vect_{\k}$:
$$H^0(Ker(\Phi(h_1)))\to H^0(\Phi(E))\to H^0(\Phi(Im(h_1))),$$
where the last arrow is $\Phi(\psi)$. Since $\Phi(\psi)(f)=0$ we conclude that there exists $f_1\in H^0(Ker(\Phi(h_1)))$ which is mapped into $f$. This contradicts to the assumption that the pair $(E,f)$ is framed stable. The Proposition is proved. $\blacksquare$

\begin{cor}
The moduli stack of stable framed objects is in fact a  scheme.
\end{cor}

In many examples it is a smooth projective scheme (cf. [Re1]).

\subsection{Stable framed objects and torsion pairs}

The above definitions can be repeated almost word by word, if we replace an ind-Artin (or locally regular) triangulated category $\CC$ by an ind-Artin abelian category ${\mathcal A}$. Then we have a definition of the framed and stable framed objects in the framework of abelian categories.
Let us discuss its relation to the classical notion of torsion pair (see e.g.  [H] for a short introduction).

Recall that a torsion pair for the abelian category ${\mathcal A}$ is given by a pairs of two full subcategories ${\mathcal T}, {\mathcal F}\subset {\mathcal A}$ such that $Hom(T,F)=0$ for any pair $T\in Ob({\mathcal T}), F\in Ob({\mathcal F})$ and such that any object $E\in Ob({\mathcal F})$ admits (a unique) decomposition
$$0\to T\to E\to F\to 0$$
with the same meaning of $F$ and $T$. Here $T$ is called the {\it torsion} part of $E$ and $F$ is called the {\it torsion-free} part of $E$. The origin of the terminology is clear from the theory of abelian groups or theory of coherent sheaves on curves.

Let us assume as before that our abelian category${\mathcal A}$  is $\k$-linear.
Suppose we are given a stability condition on ${\mathcal A}$ with the central charge $Z$. Fix $\theta\in (0, \pi)$. Then the pair of full subcategories ${\mathcal T}_{\theta}=\{T\in Ob({\mathcal A}|Arg(Z(T))>\theta\}$, ${\mathcal F}_{\theta}=\{F\in Ob({\mathcal A}|Arg(Z(F))\le  \theta\}$ defines a torsion pair for ${\mathcal A}$ (one can exchange strict and non-strict inequality signs). Let us fix a  non-zero object $P\in Ob({\mathcal A})$. It defines a functor ${\mathcal F}_{\theta}\to Vect_{\k}$ given by $\Phi(E)=Hom(P,E)$. Framed objects are pairs $(E,f:P\to E)$. Then we can give the following version of the notion of stable framed object: {\it $(E,f)$ is stable framed if either $f$ is epimorphism or $Coker(f)$ is a non-zero object of  ${\mathcal T}_{\theta}$}.

Then the above Proposition 3.1.3 still holds, and the proof is much simpler.

\begin{prp} The automorphism group of a stable framed object is trivial. 
\end{prp}
 
{\it Proof.} Let $h: E\to E$ be an automorphism such that $h\circ f=f$.  Then $(h-id)$ vanishes on the image of $f$. If $f$ is an epimorphism, we conclude that $h=id$. Otherwise, assume $h\ne id$. Then $(h-id)$ defines a non-trivial morphism $Coker(f)\to E$ which contradicts to the assumption on $Coker(f)$ and the definition of torsion pair. Hence $h=id$.
$\blacksquare$
 
 From this Proposition we again conclude that stable framed objects form a  scheme, not a stack.

\begin{rmk}
Notice that in the proof of the  Proposition 3.2.1 we did not really use a fixed slope  $\theta$, we rather worked with an individual object $E$. Hence we can give the following version of the notion of stable framed object for the framing functor defined by means of an object $P$: stable framed object is a pair $(E,f)$ such that $E$ is a non-zero object of  category ${\mathcal A}$, and $f:P\to E$ is a morphism which is either an epimorphism or a morphism with non-zero cokernel satisfying the condition that $Arg(Coker(f))>Arg(E)$ (we denote $Arg(Z(E))$ by $Arg(E)$ to simplify the notation). Yet another possibility is to require that all Harder-Narasimhan factors of $E$ belong to ${\mathcal T}_{\theta}$ (or require that all HN factors of $E$ has arguments strictly bigger than the $Arg(E)$). For all described versions the Corollary 3.1.4 remains true.
\end{rmk}

\subsection{Stable framed representations of quivers}
Let $\k$ be an algebraically closed field.

In the case of quivers without potential there is a well-known way (exploited by Nakajima and Reineke among others) to construct framed objects by adding a new vertex $i_0$ and $d_i$ new arrows $i_0\to i$ for each vertex $i\in I$ of the quiver $Q$. If we denote by $W_i$ the vector space spanned by $d_i$ arrows, then the framing functor $\Phi$ assigns to a representation $E=(E_i)_{i\in I}$ the vector space $\prod_{i\in I}Hom(W_i,E_i)$. Let $\gamma=(\gamma^i)\in \Z_{\ge 0}^I$ be a dimension vector.

Then a  framed representation of $Q$  is given by a representation  of the extended quiver $\widehat{Q}$ with the set of vertices $I\sqcup \{i_0\}$ of dimension $(\gamma^{i_0}=1,(\gamma^i)_{i\in I})$, a collection of new $d_i$ arrows $i_0\to i, i\in I$,  and a collection of linear maps $W_i\to E_i$.

Let us fix a  central charge  $Z:\Z^I\to \C$ and a ray  $l:=l_{\theta}=\R_{>0}e^{i\theta}, 0< \theta\le \pi$. Recall that we have the category $\CC_l$  of semistables with the central charge in $l$. A framed representation is  stable framed  if the following condition is satisfied (see e.g. [Re1]):

{\it  the representation $E$ of the  quiver $Q$ is semistable with central charge in $l$, and satisfies the condition that it does not have a subrepresentation $E^{\prime}$ which contains the images of all vector spaces $W_i, i\in I$ and has a bigger  argument of the central charge}.

There are many versions of the above criterion. For example, one can start with several additional vertices instead of just one. Also, one can restate the above criterion in terms of stable representations of the extended quiver $\widehat{Q}$. The later approach makes it clear why the notion of stable framed representation can be thought of as a generalization of the notion of a cyclic representation.

\begin{rmk}
For the quiver $Q_2$ with one vertex and two loops  there are no nontrivial stability conditions. Then  stable framed objects is the same as  left ideals of finite codimension in the path algebra of $Q_2$. The moduli space of stable framed objects is known as the {\it non-commutative Hilbert scheme of $\k^2$}.
\end{rmk}

\section{Modules over COHA from  stable framed objects}

\subsection{Quiver case}

Let $\k$ be an algebraically closed field.

Let fix a quiver $Q$ with the set of vertices $I$ as well as a central charge $Z: \Z^I\to \C$. We also fix  a slope $0<\theta\le \pi $ and the corresponding ray $l=l_{\theta}=\R_{>0}\cdot e^{i\theta}$. In order to specify the framing we fix a collection $(d_i)_{i\in I}$ of non-negative integer numbers. An additional (framing) vertex is denoted by $i_0$.
The corresponding extended quiver will be denoted by $\widehat{Q}:=Q^{i_0}((d_i)_{i\in I})$.

Given a dimension vector $\gamma\in \Z_{\ge 0}^I$ we denote by $\M_{\gamma,(d_i)_{i\in I}}^{st}:=\M_{\gamma,(d_i)_{i\in I}}^{st,l}$ the scheme of stable framed representations of dimension $\gamma$ having $Z(\gamma)\in l$. We denote by $\M_{\gamma,(d_i)_{i\in I}}:=\M_{\gamma,(d_i)_{i\in I}}^{l}$ the bigger space of framed representations (no stability conditions is imposed).
The group $\G_{\gamma}=\prod_iGL(\gamma^i,\k)$ acts freely on $\M_{\gamma,(d_i)_{i\in I}}^{st}$. We denote by $V_{\gamma,(d_i)_{i\in I}}^{l}=V_{\gamma,(d_i)_{i\in I}}^{\theta}$ the graded vector space $H^{\bullet}_{\G_{\gamma}}(\M_{\gamma,(d_i)_{i\in I}}^{st})=H^{\bullet}(\M_{\gamma,(d_i)_{i\in I}}^{st}/\G_{\gamma})$.

Recall that with the  ray $l=l_{\theta}$ we can associate  COHA $$\mathcal{H}_l=\oplus_{\gamma\in \Z_{\ge 0}^I, Z(\gamma)\in l}H^{\bullet}_{\G_{\gamma}}(\M_{\gamma}^{ss}).$$

Let us denote by $S:=S_{\gamma_1,\gamma_2,\gamma_3,(d_i)_{i\in I}}$ the scheme of short exact sequences
$$0\to E_1\to E_2\to E_3\to 0$$
of such representations of the extended quiver $\widehat{Q}$ that $dim(E_i)=\gamma_i\in \Z_{\ge 0}^I, i=1,2,3$, representation $E_1$ is framed, $E_2,E_3$ stable framed, and the morphism $E_2\to E_3$ is equal to the identity at the vertex $i_0$.

There is a projection $\pi_{13}:S\to \M_{\gamma_1}\times \M_{\gamma_3,(d_i)_{i\in I}}^{st}$
which sends the short exact sequence $0\to E_1\to E_2\to E_3\to 0$ to the pair $(E_1,E_3)$, where we treat $E_1$ as a representation of $Q$.   Similarly we have a projection $\pi_2$ to $E_2$. Notice that the latter is a proper morphism of $S$ to $\M_{\gamma_2,(d_i)_{i\in I}}^{st}$. Since the automorphism group of the moduli space of stable framed objects is trivial, we see that the morphism $\pi_{2\ast}\pi_{13}^{\ast}$ gives rise to a map of  cohomology groups
$$H^{\bullet}_{\G_{\gamma_1}}(\M_{\gamma_1}) \otimes H^{\bullet}( \M_{\gamma_3,(d_i)_{i\in I}}^{st})
\to H^{\bullet}(\M_{\gamma_2,(d_i)_{i\in I}}^{st}).$$

\begin{prp}
The above map gives rise to a (left) $\mathcal{H}_l$-modules structure on the vector space $V^l:=V_{(d_i)_{i\in I}}^{l}=\oplus_{\gamma}V_{\gamma,(d_i)_{i\in I}}^{l}$.
\end{prp}

{\it Proof.} Similar to the proof of associativity of the product on COHA given in [KoSo5]. $\blacksquare$

\begin{rmk} The above considerations can be generalized to the case of quivers with potential.

\end{rmk}

\begin{exa} In the case of the quiver $Q_2$ (one vertex and two loops) and $d_1=1$ the moduli space $\M_{\gamma,d_1}^{st}, \gamma\in \Z_{\ge 0}$ is the same as the moduli space of representations of the free algebra $\k\langle x_1,x_2\rangle$ of dimension $\gamma$ which are cyclic. In other words, it is the moduli space of codimension $\gamma$ ideals in the free algebra with two generators, i.e. it is the non-commutative Hilbert scheme. The above Proposition claims that it carries a structure of module over the COHA for $Q_2$ (which is the infinite Grassmann algebra). Explicit formulas for this module structure (and their generalization to the case of arbitrary number of loops) can be found in [Fra].

\end{exa}

Consider  COHA $\mathcal{H}$ of a quiver which has at least one vertex $i_0$ without loops. Then $\mathcal{H}$ is a module over the infinite Grassmann algebra (a.k.a free fermion algebra) $\Lambda^{\bullet}$. Indeed, consider $i_0$ as a quiver $Q_0$ (one vertex, no loops). We know that COHA of $Q_0$ is $\Lambda^{\bullet}$. Since it is a subalgebra of $\mathcal{H}$, it acts on $\mathcal{H}$ by left multiplication.

Let $Q$ be a quiver  with the set of vertices $I$. Let us fix a  set of non-negative integers $d=(d^i)_{i\in I}$ (not all equal to zero)and the dimension vector $\gamma=(\gamma^i)_{i\in I}$. Then we have an extended quiver $\widehat{Q}$ with the set of vertices $I\sqcup i_0$ and $d^i$ arrows from $i_0$ to $i\in I$. For a fixed central charge $Z: \Z^I\to \C$ the moduli space $\M_{\gamma,d}^{st,l}$ of  stable framed representations of $Q$ of dimension $\gamma$ such that $Z(\gamma)\in l$ is a non-empty smooth  variety of pure dimension $\sum_{i\in I}d^i\gamma^i-\chi(d,d)$, where $\chi(\alpha,\beta)$ is the Euler-Ringel bilinear form of $Q$ (see [EnRe], Prop. 3.6).
Moreover it admits a projective morphism to the moduli space of polystables of fixed slope.

\subsection{Representations of COHA in general case}
We will give a sketch of the construction.

Let ${\mathcal A}$ be ``good' abelian subcategory in the $3CY$ category ${\mathcal C}$. We assume the conditions on the potential $W$ which guarantee existence of COHA of ${\mathcal A}$ as well as moduli spaces of stable framed objects.   Then considerations from the previous subsection can be generalized to this situation provided ${\mathcal A}$ satisfies some extra conditions, e.g. that classes $cl(E)$ of objects of ${\mathcal A}$ belong to an additive monoid $\Gamma_+$ which is mapped to $\Z_{\ge 0}^n$ under the chosen identification $\Gamma\simeq \Z^n$.

Next, let us fix a ray $l=\R_{\ge 0}\cdot e^{i\theta}$ in the upper half-plane, and a stability function $Z:\Gamma\to \C$ such that $Z(\Gamma_+)$ belongs to the upper half-plane. Then we have the category ${\aA}_l$ of semistables with the central charge in $l$. Let us fix the framing functor $\Phi$.
Then we can speak about framed and stable framed objects.

Recall that there is a notion of morphism of framed objects $(E_2,f_2)\to (E_3,f_3)$.  An epimorphism $(E_2,f_2)\to (E_3,f_3)$ is a morphism in the category of framed objects which induces a homomorphism $H^0(\Phi(E_2))\to H^0(\Phi(E_3))$ which sends $f_2$ to $f_3$ (see Section 3.1 for the notation).

Assume that $E_2$ and $E_3$ are semistable objects with central charges in the ray $l$. Then the kernel of the epimorphism $(E_2,f_2)\to (E_3,f_3)$ of framed objects does not have to be framed. Let us consider the stack ${\mathcal Z}_{\gamma_1,\gamma_2}$ of triples $(E_1,(E_2,f_2),(E_3,f_3))$ where:

a) $cl(E_1)=\gamma_1$, and $Z(\gamma_1)\in l$;

b) $(E_2,f_2)$ is stable framed, $cl(E_2)={\gamma_1+\gamma_2}$, $Z(cl(E_2)\in l$;

c) $(E_3,f_3)$ is stable framed, $cl(E_3)=\gamma_2$, $Z(cl(E_2)\in l$;

d) there is a epimorphism of framed objects $(E_2,f_2)\to (E_3,f_3)$ such that it induces (in the category of semistable objects with the central charge in $l$) a short exact sequence
$$0\to E_1\to E_2\to E_3\to 0.$$

Recall that stable framed objects with fixed class $\gamma\in \Gamma_+$ form a scheme which we denote by $\MM_{\gamma}^{st,fr}$. Then we have natural projections $p_2: {\mathcal Z}_{\gamma_1,\gamma_2} \to \MM_{\gamma_1+\gamma_2}^{st,fr}$ and $p_3: {\mathcal Z}_{\gamma_1,\gamma_2} \to \MM_{\gamma_2}^{st,fr}$ which are morphisms of stacks.
Furthermore, let $\MM_{\gamma}$ denotes the moduli stack of objects of ${\aA}_l$. Then we have the natural projection $p_1: {\mathcal Z}_{\gamma_1,\gamma_2} \to \MM_{\gamma_1}$. 

We will assume that:

i) if we consider the analog of the above situation with all $f_i=0, i=1,2,3$ (i.e. we work just in the abelian category $\aA_l$) then the restriction of $p_2$ to $p_1^{-1}(E_1)\cap p_3^{-1}(E_3)$ is a morphism of smooth proper stacks;

ii) in general, for fixed $E_1$ and $(E_3,f_3)$ as above, the restriction of $p_2$ to $p_1^{-1}(E_1)\cap p_3^{-1}((E_3,f_3))$ is a morphism of smooth proper stacks.

By condition i) COHA ${\mathcal H}_l$ of the category $\aA_l$ is well-defined as an associative algebra.
For that we use the critical version of the cohomology from [KoSo5] with trivial potential.
Furthermore, repeating the construction from the previous subsection we obtain  a structure of (left) ${\mathcal H}_l$-module on $V:=V^l=\oplus_{\gamma\in \Gamma_+}H^{\bullet}(\MM_{\gamma}^{st,fr})$.

\begin{rmk}

More generally, we can construct modules over  COHA  by considering the stack of objects whose Harder-Narasimhan filtration has consecutive factors with arguments of the central charge  belonging to the interval $[\theta, \pi]$. If we have an exact short sequence
$$0\to E_1\to E_2\to E_3\to 0,$$
such that  $Arg\,Z(E_2)\in [\theta,\pi]$ then $Arg\,Z(E_3)$ belongs to the same interval, while $Arg\,Z(E_a)\in [0,\pi]$. Then we get a representation of COHA in the cohomology of the stack of objects generated by semistables
$E$ such that $Arg\,Z(E)\in [\theta,\pi]$.

Similarly, one can show that if $V$ is a strict sector in the plane then the graded ``Cohomological Hall vector space" $\mathcal{H}(V)$ bounded from the left by a ray $l$ is a module over the COHA $\mathcal{H}_l$ associated with the ray.

\end{rmk}

Furthermore, suppose that our abelian category ${\aA}$ is a ``good" subcategory of an ind-Artin $3CY$ category $\CC$ endowed with orientation data. Let $W$ be the potential for $\CC$. It gives rise to the sheaf of vanishing cycles $\phi_W$ on the stack of objects of $\CC$. Then the pull-backs of $\phi_W$ to the stack of objects of $\aA$ and subsequently to  $\MM_{\gamma}$ and $\MM_{\gamma}^{st,fr}$ are well-defined. Then, similarly to [KoSo5] (and under the above assumption), the above construction (but this time with cohomology groups with coefficients in $\phi_W$) gives rise to the module $V:=V^l=\oplus_{\gamma\in \Gamma_+}H^{\bullet}(\MM_{\gamma}^{st,fr},\phi_W)$ over the COHA ${\mathcal H}_l$ of $\aA_l$. The details will be explained elsewhere.

\subsection{Hecke operators associated with simple objects}

In the classical Nakajima construction of the infinite Heisenberg algebra (see [Nak2]) one considers pairs of ideal sheaves $(J_2,J_3)$ on a surface $S$ such that $J_2\subset J_3$ and $Supp(J_3/J_2)=\{x\}$, where $x$ is a fixed point.
Then one has an epimorphism ${\mathcal O}_S/J_2\to {\mathcal O}_S/J_3$.
Let us compare this observation with the above construction of modules over COHA. We see that a fixation of $K$-theory classes $\gamma_i, i=1,2$ for a pair of objects $(E_2,E_3)$ along with an epimorphism $E_2\to E_3$ corresponds in the Nakajima's construction to the fixation of $n_i, i=2,3$ such that $J_i\in Hilb_{n_i}(S)$ and to the above-mentioned epimorphism of the quotient sheaves.  

In the construction  of the module structure on the cohomology of stable framed objects we used the pushforward map associated with the projection to the middle term in the moduli space of short exact sequences
$$0\to E_1\to E_2\to E_3\to 0,$$
where $E_2,E_3$ are stable framed (we omit here $f_i, i=2,3$ from the notation).
As a result, our construction gives rise to the ``raising degree" operators ${\mathcal H}_{\gamma_1}\otimes V_{\gamma_2}\to V_{\gamma_1+\gamma_2}$ for the COHA action ${\mathcal H}\otimes V\to V$. There are no ``lowering degree" operators, which  would correspond to the projection to the term $E_3=E_2/E_1$. The reason is similar to the one in the Nakajima's construction: such a projection is not proper.

Originally Nakajima solved the problem by considering points $x\in S$ which belong to a compact subset in $S$. We can use this idea and consider short exact sequences as above, where $E_1$ is a {\it simple} object which runs through a compact (in analytic topology) subset in the moduli scheme of simple objects of our abelian category ${\mathcal A}$.

Let us illustrate the construction in the case of quivers without potential and trivial stability condition. In that case stable framed objects are  cyclic modules over the path algebra of the quiver.  Then we should prove that there are sufficiently many cyclic modules with the fixed simple submodule and fixed cyclic quotient. This is guaranteed by the following result.

\begin{prp} Let $A$ be an associative algebra, $(M_2,v_2),(M_3,v_3)$ be   $A$-modules with marked elements $v_i\in M_i, i=2,3$ such that $v_3$ is a cyclic vector for $M_3$. Let $f:M_2\to M_3$ be an epimorphism of $A$-modules such that $f(v_2)=v_3$ and such that $W=Ker(f)$ is a simple $A$-module. Suppose that the extension
$$0\to W\to M_2\to M_3\to 0$$
is non-trivial.
Then $v_2$ is a cyclic vector for $M_2$.

\end{prp}
{\it Proof.} Let $M_2^{\prime}\subset M_2$ be the $A$-submodule generated by $v_2$. If $M_2^{\prime}=M_2$ then we are done. Otherwise we have a non-trivial epimorphism $g:W\to M_2/M_2^{\prime}$ of $A$-modules. Its kernel is a submodule of $W$. It must be trivial, since $W$ is simple. Hence $g$ is an isomorphism. Then the submodules $W$ and $M_2^{\prime}$ determine the direct sum decomposition $M_2=W\oplus M_2^{\prime}$, where $M_2^{\prime}\simeq M_3$. Hence the extension
$0\to W\to M_2\to M_3\to 0$ is trivial. This contradiction shows that $v_2$ is a cyclic vector. $\blacksquare$

\begin{cor}
For fixed $W,M_3$ the stack of cyclic modules $M_3$ which are middle terms in the above short exact sequence is  a smooth projective scheme isomorphic to the projective space ${\bf P}(Ext^{1}(M_3,W))$.

\end{cor}

{\it Proof.} Follows from the Proposition. $\blacksquare$

\begin{rmk}
Similar result holds in case when $M_3$ is stable framed and $S$ is simple. 
\end{rmk}

Let now $\MM^{simp}:=\MM^{simp}_{\mathcal A}$ be the moduli space of simple objects in the heart ${\cal A}$ of the ``good" $t$-structure of an ind-Artin $3CY$ category $\CC$ endowed with orientation data. Then $\MM^{simp}_{\mathcal A}$  is a smooth separated scheme. Let $H^{\bullet}_{BM,c}(\MM^{simp})$ denotes  compactly supported Borel-Moore cohomology. As before we have two projections $\pi_1,\pi_3$ from the  schemes of short exact sequences to its first and last term, i.e. to the moduli space  $\MM^{simp}$ of simple objects and to the moduli space  $\MM^{st}$ of stable framed objects correspondingly. Then the composition $\pi_{3,\ast}\circ \pi_1^{\ast}$ defines a collection of operations on
$H^{\bullet}(\MM^{st})$ parametrized by the elements of $H^{\bullet}_{BM,c}(\MM^{simp})$. The above Proposition (or rather its analog for non-trivial stability condition) ensures that the operations are well-defined. Differently from the action of COHA defined in the previous subsection, these operations {\it decrease} the degree $\gamma\in \Gamma$.

\begin{rmk}
Let us recall that for any $i\in \pm \Z_{>0}$ Nakajima defines an operator $P[i]$ which corresponds to the $i$-th generator of the infinite Heisenberg algebra. In the above discussion the operator $P[i]$  corresponds to the direct sum $iS:=S\oplus S\oplus ....\oplus S$ of $\pm i>0$ of copies of the simple object $S$.
\end{rmk}

Using the above construction one can extend a representation of COHA to a representation of a bigger algebra, which we call ``full COHA" (or double of COHA). We do not know yet how to define this algebra intrinsically. 
Our approach is similar to the Nakajima's construction of the infinite Heisenberg algebra from two representations of the symmetric algebra: one is given by creation operators and another one is given by annihilation operators. Commuting creation and annihilation representations {\it in the representation space} Nakajima recovers the infinite Heisenberg algebra.
One can also compare the above construction with the one in [Re1]. 

\begin{rmk}
As we already mentioned, differently from the the case of constructible Hall algebras, we do not know a compatible comultiplication on COHA. This would help to define full COHA by means of the Drinfeld double construction. Having in mind that in the case of quivers without potential COHA is a shuffle algebra, one can hope for explicit formulas similar to those in [Neg].
\end{rmk}

\subsection{Full COHA: an example}

It is well-known that one can obtain finite-dimensional representations of quantized enveloping algebra
of finite-dimensional semisimple Lie algebras using framed stable representations of quivers and constructible Hall algebra  (see e.g. [Re1] where this idea which goes back to Nakajima was implemented for  representations of the ``positive" part of a quantum group).

We are going to  recall below how to recover quantized enveloping algebra $U_q(sl(2))$ by combining two representations of the Hall algebra for quiver $A_1$.

Recall that if for a quiver $Q$ we take the stability function $\Theta=0$, then every finite-dimensional representation of $Q$ is semistable. The  moduli space of stable framed representations admits in this case a simple description in terms of Grassmannians (see e.g. [Re1], Prop. 3.9).

For the quiver $A_1$ (one vertex $i_1$ and no arrows) the framing consists of  a new vertex $i_0$ and $d$ arrows $i_0\to i_1$. The stability function is trivial automatically, and one can easily see that for each dimension vector $\gamma\in \Z_{\ge 0}$ the moduli space ${\cal M}_{\gamma,d}=\M_{\gamma}^{\theta=0,st}$ of framed stable representations of dimension $\gamma$ is isomorphic to the  Grassmannian $Gr(d-\gamma,d)\simeq Gr(\gamma,d)$. Hence it is non-empty  for $\gamma\le d$ only. Let us denote by $Gr(d)$ the ``full Grassmannian" consisting of vector subspaces of $\C^d$ of all dimensions (this space is disconnected).
Then the moduli space of $d$-framed semistable representations of $A_1$ is $Gr(d)$.

The space of $GL(d)$-invariant functions with finite support $Fun^{GL(d)}(Gr(d))$ is a module over the constructible Hall algebra of $A_1$. 
The constructible Hall algebra of the quiver $A_1$ is the polynomial algebra with one generator 
$z:={\bf 1}_{\C}$, where the generator $z$ corresponds to the characteristic function ${\bf 1}_1$ of $\M_{1}$ in the stack $\MM=\sqcup_{\gamma\ge 0}\M_{\gamma}$. Indeed the Hall product gives an isomorphism of the constructible Hall algebra with the polynomial ring $\C[z]$. 
In each $Gr(k,d)$ we have only one $GL(d)$-orbit of the standard coordinate vector subspace $\C^k\subset \C^d$. Let us denote by $v_k, 0\le k\le d$ the characteristic function of the corresponding $GL(d)$-orbit.

Let us consider the ``minus" Hecke correspondence given by pairs $(V_{k-1}\subset V_k)$ with $1$-dimensional factor $V_1$ and project to $V_{k-1}$. Equivalently, we consider the projection to the first terms from  the set of short exact sequences
$$0\to V_{k-1}\to V_{k}\to V_{1}\to 0.$$
Using the pullback/pushforward construction discussed previously,
we obtain a representation of the algebra $\C[z]$ given by $\rho_{-}(z)v_k={{q^{k}-q^{-k}}\over {q-q^{-1}}}v_{k-1}, 1\le k\le d$, and $\rho_{-}(z)v_0=0$, where the factor  comes from the normalization of the cocycle $c(M,N)$ above as $q^{\chi(M,N)}$. The Euler-Ringel form $\chi$  on the pair of representations $E$ of dimension $a$ and $F$ of dimension $b$ is given by $\chi(E,F)=ab$.
Similarly, consider the ``plus" Hecke correspondence $(V_k\subset V_{k+1})$ and project to $V_{k+1}$. Then we get a representation of $\C[z]$ in ${\cal F}_n$ given by
$$\rho_+(z)v_k={{q^{k+1}-q^{-k-1}}\over{q-q^{-1}}}v_{k+1}, 0\le k\le d-1, \rho_+(z)v_d=0.$$
Combining $\rho_{-}$ and $\rho_{+}$ together we obtain the standard $d$-dimensional representation of the quantized enveloping algebra $U_q(sl(2))$ where the ``positive" generator $E$ is represented by $\rho_{-}(z)$ while the ``negative" generator $F$ is represented by $\rho_{+}(z)$. Then the commutator $[E,F]$ maps $v_k$ to ${q^{2k}-q^{2k}\over{q-q^{-1}}}v_k$. From this formula one can recover the action of the Cartan generators $K,K^{-1}$.

Let us apply similar considerations in the case of COHA of the quiver $A_1$. Recall, in this case COHA is isomorphic to the algebra $\Lambda^{\bullet}=\Lambda^{\bullet}(\xi_1,\xi_2,....), deg\,\xi_{i}=2i-1,i\ge 1$. Since for the trivial stability function the COHA associated with the ray $\theta=0$ coincides with whole COHA, we obtain a representation of the infinite Grassmann algebra $\Lambda^{\bullet}$ in the finite-dimensional vector space $V:=H^{\bullet}(Gr(d))=\oplus_{0\le k\le d}H^{\bullet}(Gr(k,d))$. One can write down explicitly the action of the  generators of $\Lambda^{\bullet}$ on the cohomology classes of Schubert cells. 
 
Following the general definition, we consider the moduli stack of short exact sequences 
$$0\to E_1\to E_2\to E_3\to 0,$$
where $E_2$ and $E_3$ are stable framed of the same slope, and $E_1$ is a representation without framing.  First we consider the representation of COHA coming from the projection to $E_2$.

In order to write down the corresponding representation explicitly let us choose a subspace in each $H^{\bullet}(Gr(k,d)), 0\le k\le d$ spanned by the cohomology classes corresponding to $(\C^{\ast})^d$-fixed points. We denote this basis by ${e}_{j}:={\bf 1}_{\C_{j_1,...j_k}}$ (recall that the fixed points correspond to coordinate subspaces $\C_{j_1,...j_k}\subset\C^d$ spanned by the standard basis vectors $f_{j_1},...,f_{j_k}, j_1<j_2<...<j_k$). We can identify the graded vector space $V$ with the quotient $\Lambda^{\bullet}(\xi_1,...,\xi_d)/I_d$, where $I_d$ is the ideal generated by  $\xi_{i}, i\ge d+1$. Then the pullback/pushforward construction gives us the representation of $\Lambda^{\bullet}$ in $V$ by ``creation" operators: $a_n^{\ast}:{e}_{j}\mapsto \xi_n\wedge e_{j}$.

Second, we consider the representation of COHA coming from the projection to $E_3$.
This gives a representation of $\Lambda^{\bullet}$ on $V$  by ``annihilation" operators $a_n:e_j\mapsto i_{\xi_n}(e_j)$, where $i_{\xi_n}$ is the contraction operator which delete the variable $\xi_n$ from the monomial $e_j$.

Then,similarly to the consideration with the constructible Hall algebra, we combine both representations of COHA into a single one. In this way recovers the representation  of the Lie algebras $D_{d+1}$.\footnote{I thank to Xinli Xiao for making explicit computations.} 
This leads to the following conjecture.

\begin{conj}
Full COHA for the quiver which has one vertex and $m$ arrows is isomorphic for even $m$ to  the infinite Clifford algebra $Cl_c$ with generators $\xi_n^{\pm}, n\in 2\Z+1$ and the central element $c$, subject to the  anticommuting relations between $\xi_n^+$ (resp. $\xi_n^-$) as well as the relation $\xi_n^+\xi_m^-+\xi_m^-\xi_n^+=\delta_{nm}c$. 

In the case of odd $m$ full COHA is isomorphic to the infinite (graded) Heisenberg algebra (in the above formulas change anticommuting brackets  by the commuting ones).

\end{conj}

In the case of finite-dimensional representations we have $c\mapsto 0$ and  two representations of the infinite Grassmann algebra, which are combined in the representations of the orthogonal Lie algebra as explained above. Notice that $Cl_c$ is the Clifford algebra associated with the positive part of the affine Lie algebra $sl(2)$. This might explain  the relation of the full COHA to the quiver $A_1$.

\section{Some representations of COHA motivated by physics and knot theory}
In this section we are going describe some interesting classes of representations of COHA. Details of the constructions will appear elsewhere. The reader can consider this section as a collection of speculations. 

\subsection{Fukaya categories of conic bundles and gauge theories from class ${\mathcal S}$}

We are going to illustrate the idea in the case of   $SL(2)$ Hitchin integrable systems. Our motivation is the general conjecture (F.1) from the Introduction of  [ChDiManMoSo]. In this particular case it admits a very precise interpretation. Namely, with a point of the universal cover of the  base of Hitchin system on a complex curve, say, $C$ one can associate a compact Fukaya category endowed with a stability structure (``compact'' means that it is generated by local systems supported on compact Lagrangian submanifolds). It is the Fukaya category of a non-compact Calabi-Yau $3$-fold $X$ described by the corresponding spectral curve (see [KoSo8] for a more general framework). The compact Fukaya category is endowed with the natural t-structure generated by SLAGs which are 3-dimensional Lagrangian spheres. The central charge of the corresponding stability condition is given by the period map of the Liouville form on $T^{\ast}C$ restricted to the spectral curve. According to the general theory developed in Section 8 of [KoSo1] categories generated by spherical collections are in one-to-one correspondence with pairs $(Q,W)$, i.e. quivers with potential. Hence we can speak about corresponding COHA and its representations in the cohomology of the moduli spaces of stable framed objects of the category $Crit(W)$. This would give an interesting class of representations of the BPS algebra of the corresponding gauge theory from class ${\mathcal S}$.

Recall that surface defects in physics correspond to points of the curve $C$. In terms of the corresponding non-compact Calabi-Yau $3$-fold they are complex $2$-dimensional submanifolds of $X$. Consider the the moduli space of SLAGs with the boundary which belongs to such a submanifold $S$. The corresponding category ${\mathcal F}(X,S)$ can be thought of as a version of Fukaya-Seidel category of thimbles (see [Se]) with the (analog of the) potential being the natural map $X\to C$.

Furthermore, the operation of connected Lagrangian sum plays a role of an extension in the compact Fukaya category ${\mathcal F}(X)$. This operation underlies the product structure on the COHA ${\mathcal H}^{Q,W}$. 

Let us observe that there is an operation of taking the connected Lagrangian sum of a Lagrangian submanifold without boundary and the one with the boundary on $S$. Mimicking the definition of the product on COHA with the ``moduli space of Lagrangian connected sums'' instead of the subvariety $\M_{\gamma_1,\gamma_2}\subset \M_{\gamma_1+\gamma_2}$ (see [KoSo2], Section 2), one obtains the ${\mathcal H}^{Q,W}$-module structure on the cohomology of the moduli space of SLAGs with the boundary on $S$. \footnote{Alternatively, following Paul Seidel, one can consider the double cover of the Calabi-Yau $3$-fold branched along  the divisor given by the complex surface. Then Lagrangian submanifolds with boundary lift to closed ones in the branched 
cover. One can form an equivariant Lagrangian connected sum, and then
interpret it as an operation 
on the original Lagrangian submanifolds with boundary.}

\subsection{Resolved conifold and quivers}

Let $X=tot({\cal O}(-1)\oplus {\cal O}(-1))$ be the resolved conifold. We denote the zero section of the corresponding vector bundle by $C_0\simeq {\bf P}^1$. Let  us fix a point $p_0:=0\in C_0$.

Let ${\mathcal A}$ be the abelian category  of {\it perverse coherent sheaves} on $X$ topologically supported on $C_0$ (see e.g. [NagNak], [Tod] for descriptions convenient for our purposes; in [NagNak] our category ${\mathcal A}$ was denoted by $Per_c(X/Y)$, where $X\to Y$ is the crepant resolution of the conifold singularity  $Y=\{xy-zw=0\}$).

It is known (see e.g. [NagNak]) that ${\mathcal A}$ is equivalent to the abelian category  $Crit(W)$ associated with the pair $(Q,W)$, where $Q$ is a quiver with two vertices $i_1,i_2$ two arrows $a_1,a_2:i_1\to i_2$, two arrows $b_1,b_2:i_2\to i_1$ and ``Klebanov-Witten potential" $W=a_1b_1a_2b_2-a_1b_2a_2b_1$. In particular, for any $\gamma=(\gamma^1,\gamma^2)\in \Z_{\ge 0}^2$ the stack of objects of $Crit(W)$ of dimension $\gamma$ is equivalent to the stack of such representations of $Q$ of dimension $\gamma$ in coordinate vector spaces, which belong to the critical locus of the function $Tr(W)$.

We recall that the category of perverse coherent sheaves carries a family of  geometrically defined {\it weak stability conditions} (see e.g. [Tod]).   In the case of the category $Crit(W)$ there is a class of  stability conditions associated with the slope function.

Equivalence of these  two categories gives rise to the ``chamber" structure of the space of stability conditions on ${\mathcal A}$ described in [NagNak]: some of the (infinitely many) chambers correspond to the quiver-type stability conditions, while ``at infinity" we have chambers of geometric origin corresponding to different choices of the weak stability condition.  

One has a similar story when the framing is taken into account. Then one deals with framed perverse coherent sheaves and framed representations of $(Q,W)$ (i.e. critical points of the function $Tr(W)$ considered as a function on the space of representations of the extended 
quiver $\widehat{Q}$ obtained from $Q$ by adding an extra vertex $i_0$ and an arrow $i_0\to i_1$.)

In the ``quiver chamber", we can (after a choice of a stability condition on ${\mathcal A}$ which belongs to the above class) speak about the moduli space of stable framed objects of the fixed slope. For a given $(l,n)\in \Z_{\ge 0}\times \Z$ and a choice of  certain stability condition on the category $Rep (Q)$ of finite-dimensional representations of $Q$, and a choice of an angle $\theta$ (which depends on $(l,n)$), the space of stable framed representations of $(Q,W)$ with the slope $\theta$ becomes isomorphic to the moduli space of Pandharipande-Thomas stable pairs $P(l,n)$ (see [NagNak]). There is no single $\theta$ which serves all $(l,n)$. 

It follows from the previous section that:

\begin{prp} For a choice of stability conditions in the ``quiver chamber", COHA ${\mathcal H}^{(Q,W)}$ acts on the cohomology of the moduli space of stable framed representations of $(Q,W)$ having fixed slope.

\end{prp}

Let us observe that if we have a morphism $f: E_2\to E_3$ of PT stable pairs which is surjective in degree zero (i.e. on the sheaves supported on $C_0$) then $Ker(f)$ is a coherent sheaf scheme-theoretically supported on $C_0$. Passing to the cohomology groups we reformulate the above Proposition by saying that  COHA of the corresponding category acts on the cohomology of the moduli space of PT stable pairs.
We expect the same result to hold in ``geometric chambers", were one uses weak stability conditions.

\subsection{Vertically framed sheaves and algebraic knots}

Motivated by [GuSchVa] and many other physics papers on the relation between knot theory and BPS states for the resolved conifold, one can hope for an application of the representation theory of COHA of the resolved conifold to knot invariants. Among mathematical motivations we can mention  the main conjecture from [ORS], its reformulation in [DiHuSo] and its proof  in [Mau1] (in the ``unrefined" form). In this subsection we discuss appropriate moduli spaces following [DiHuSo] and speculate about corresponding representation of COHA.

Let $K$ be an algebraic knot or link which is obtained by Milnor construction, i.e. via the  intersection of a plane singular curve $C_K$ with the $S^3$-boundary of a small ball around the singularity.
If we would like to incorporate algebraic knots in the story,  we should add to the story coherent sheaves on $X$ supported on the  curve $C_K$ placed in the fiber of the projection  $X\to C_0={\bf P}^1$.

 More precisely we consider coherent sheaves which are ``vertically framed" along $C_K$ (see details [DiHuSo]). Stable vertically framed coherent sheaves provide a natural generalization of PT stable pairs from [PT1] (see also [DiShVa] for some physics arguments).

Let us  recall some general definitions and details following [DiHuSo].

Let $X=tot({\cal O}(-1)\oplus {\cal O}(-1))$ be the resolved conifold, $C$ a planar complex algebraic curve with the only singular point $p$.
%
%

The abelian category of $C$-framed perverse coherent sheaves  is a full subcategory ${\cal A}^{C}\subset D^b(Coh({X}))$. Roughly speaking, ${\cal A}^{C}$ consists of complexes $E$ of coherent sheaves on ${X}$ such that the cohomology sheaves $H^i(E)$ are non-trivial for $i\in \{0,-1\}$ only, and those cohomology sheaves are topologically supported one the union $C\cup C_0={\bf P}^1$ (see loc.cit. Section 2.2 and below for  more precise description). The category ${\cal A}^C$ is closed under extensions, and   it is a full subcategory of the category of perverse coherent sheaves ${\cal A}\subset D^b(X)$. After fixing  K\"ahler class $\omega$ on the compactification $\overline{X}$ defined in [DiHuSo], one defines a family  of weak stability conditions on ${\cal A}^C$ associated with an explicitly given family of slope functions $\mu_b:=\mu_{(\omega, b\omega)}$ described in the loc.cit.
Then one can speak about $C$-framed (semi)stable sheaves, meaning weakly (semi)stable objects of ${\cal A}^C$ with respect to the slope function $\mu_b$. \footnote{In [DiHuSo] the authors considered stable vertically framed sheaves on the compactification $\overline{X}$. The corresponding moduli spaces were projective. Considerations with non-compact submanifold $X$ gives rise to quasi-projective moduli spaces. We ignore these technicalities here.}

For ``very negative" value of $b$ the moduli space
${\mathcal P}^s_b({{X}} ,C,r,n)$ of $C$-framed $\mu_b$-stable objects $E$ with $ch(E)=(-1,0,[C]+r[C_0],n)$ is
isomorphic to the moduli space of stable framed pairs on ${X}$ in the sense of Pandharipande and Thomas which are $C$-framed.
If we move the value of $b$ from $b=-\infty$ to a small positive number (which depend on $r$) the above moduli space of $\mu_b$-stable objects experiences   finitely many wall-crossings. One of the main results of [DiHuSo] is a theorem which relates the moduli space of $\mu_b$-stable objects of ${\cal A}^C$ for small $b>0$ with the punctual Hilbert schemes from [ORS]. This relates the DT-invariants of the category of $C$-framed stable sheaves with  HOMFLY polynomials of algebraic knots. The moduli space  ${\mathcal P}^{ss}_b({{X}} ,C,r,n)$ of $\mu_b$ semistable objects is a $\C^{\ast}$-gerbe over ${\mathcal P}^s_b({{X}} ,C,r,n)$.

Let us fix $(r,n)\in \Z_{\ge 0}\times \Z$ and consider the full subcategory ${\mathcal A}^C_{r,n}\subset {\cal A}^C$ consisting of objects $E$ such that $ch(E)=(-1,0,[C]+r[C_0],n)$. Let $E_1$ be a pure dimension one sheaf on ${X}$ supported on $C_0$ (hence it belongs to ${\mathcal A}^C$ as well),  and let $E_3\in Ob({\mathcal A}^C_{r_3,n_3})$. Then  we see that the middle term $E_2$ of an extension in ${\cal A}^C$
$$0\to E_1\to E_2\to E_3\to 0$$
belongs to ${\mathcal A}^C_{r_2,n_2}$ for some $r_2,n_2$.

Let ${\mathcal M}^C:=\cup_{r,n}{\mathcal M}_{r,n}^C$ be the moduli space (stack) of the objects $E$ which belong the category ${\mathcal A}^C_{r,n}$ for some $r,n$. Let ${\cal M}^{C_0}$ be the moduli space (stack) of objects of the category $Coh_{C_0}({X})$ of coherent sheaves on ${X}$ supported on the zero section $C_0={\bf P}^1$. Let ${\mathcal N}$ be the moduli space (stack) of short exact sequences as above.

We have the following projections: $\pi_{13}: {\mathcal N}\to {\mathcal M}^C\times{\cal M}^{c_0},(E_1,E_2,E_3)\mapsto (E_1,E_3) $ and  $\pi_2: {\mathcal N}\to {\mathcal M}^C, (E_1,E_2,E_3)\mapsto E_2$.

Then we can apply the same procedure as for framed representations of quivers using the composition $\pi_{2\ast}\pi_{13} ^{\ast}$. Then e.g. in the case of COHA it gives us the module structure $H^{\bullet}({\cal M}^{C_0})\otimes H^{\bullet}({\mathcal M}^C)\to H^{\bullet}({\mathcal M}^C)$ over the COHA of the category $Coh_{C_0}({X})$, where by $H^{\bullet}$ we denote an appropriate  cohomology theory.

Now we can use the weak stability condition defined by the slope function $\mu_b$.
More precisely, let us choose a stability parameter $b$ satisfying the condition $(3.1)$ of Lemma 3.1 from [DiHuSo] and consider $\mu_b$-semistable objects $E$ of ${\cal A}^C$ such that $ch(E)=(-1,0,[C]+r[C_0],n)$, where $(r,n)\in \Z_{\ge 0}\times \Z$ is fixed. Then we can repeat the above definition but this time in the exact sequence
$$0\to E_1\to E_2\to E_3\to 0$$
we will assume that $E_2$ and $E_3$ are weakly semistable objects with respect to $\mu_b$, and $E_1$, as before, is an arbitrary coherent sheaf on ${X}$ supported on $C_0$.

There is an explicit description of $\mu_b$-semistable and $\mu_b$-stable objects of ${\cal A}^C_{r,n}$ for sufficiently small positive $b$ given in [DiHuSo], Section 3. For example a $\mu_b$-stable object $E$ fits into an exact short sequence
$$0\to E_C\to E\to {\mathcal O}_{C_0}(-1)^r\to 0,$$
where $E_C=({\mathcal O}_{{X}}\to F_C)$ is a stable pair on ${X}$ in the sense of Pandharipande and Thomas, with the sheaf $F_C$ scheme theoretically supported on $C$ (and satisfying some non-degeneracy conditions, see [DiHuSo], Proposition 3.3 for the details). Similarly, any $\mu_b$-semistable object fits into an exact sequence where instead of ${\mathcal O}_{C_0}(-1)^r$ one has a sheaf $G$ topologically supported on $C_0$ (and $ch_2(G)=r[C_0]$) which is a direct image (under the embedding $i:C_0\to {X}$) of the vector bundle
$\oplus_{1\le j\le m} {\cal O}(a_j)^{r_j}$ with $a_1>...>a_m\ge -1$. The Harder-Narasimhan filtration of $G$ (with respect to the $\omega$-slope defined by $\chi(G)/r$) therefore have consecutive factors with slopes $a_j/r$. 

Based on the above considerations one can hope  that  $C$-framed stable sheaves  play a role similar to the one played by  stable framed objects in the abelian categories. In particular, cohomology groups of the moduli spaces of $C$-framed stable sheaves  should give rise to representations of COHA of $X$. It is not clear at this time how far this idea can be developed. In fact computations made by E. Diaconescu, show that if in the short exact sequence $0\to F\to E_1\to E_2\to 0$ the terms $E_1,E_2$ are $C$-framed stable then $F$ is isomorphic to ${\mathcal O}(-2)^n$. It seems plausible that in order to obtain interesting representations of the full COHA, one should also include in considerations short exact sequences of the type $0\to E_1\to E_2\to F\to 0$, where $E_1,E_2$ are $C$-framed stable. This should lead to the representation of the full COHA in the way discussed previously. We expect that in this way we obtain affine $sl(2)$.

\begin{rmk}
The above story with $C$-framed sheaves is related to algebraic knots. As for more general knots, one can hope that the 
following picture can be made mathematically precise.

For any non-compact real analytic Lagrangian submanifold $L\subset X$ with ``good behavior at infinity" there should be a well-defined stack $Coh_{\le 1}(X,L)$ of real analytic sheaves on $X$ (considered as a real analytic manifold) with the following properties:

a) Every $F\in Coh_{\le 1}(X,L)$ has topological support, which is  an immersed $2$-dimensional real-analytic submanifold of $X$. Moreover, the support without boundary is an immersed non-compact complex analytic curve. The restriction of $F$ to the complement of the boundary is a a coherent sheaf on the corresponding complex manifold.

b) The boundary of the support of each $F\in Coh_{\le 1}(X,L)$ belongs to $L$.

c) The stack $Coh_{\le 1}(X,L)$ is a countable union of  real-analytic stacks of finite type. It is naturally the stack of objects of the abelian category of real-analytic sheaves on $X$ satisfying conditions a) and b).

In particular, sheaves $F$ with pure support are those for which the support is an immersed ``bordered Riemann surfaces" in the sense of [KatzLiu]. 

One can hope that despite of the analytic nature of objects, there is  a theory of stability structures for this category, as well as the notion of stable framed object. 

Notice that we can consider extensions $0\to F\to E_1\to E_2\to 0$, where $E_1,E_2$ are objects of $Coh_{\le 1}(X,L)$, while $F$ is the usual coherent sheaf on $X$ with support on $C_0={\bf P}^1$. We expect that this operation leads to the action of COHA on the cohomology of framed stable objects in $Coh_{\le 1}(X,L)$, similarly to the case of $C$-framed stable sheaves.

\end{rmk}

Finally, if the above discussion about representations of  COHA of the resolved conifold makes sense, then one can hope that it is related to  the representation theory of DAHA discussed in [GorORS].

\vspace{3mm}

{\bf References}

\vspace{2mm}

[BLM] A. Beilinson, G. Lusztig, R. Macpherson, A geometric setting for the quantum deformation of $GL_n$, Duke Math. J.,61:2, 1990.

\vspace{2mm}

[BuJoMe] V. Bussi, D. Joyce, S. Meinhardt, On motivic vanishing cycles of critical loci,  arXiv:1305.6428.

\vspace{2mm}

[ChDiManMoSo] W. Chuang, D. Diaconescu, J. Manschot, G. Moore, Y. Soibelman, Geometric engineering of (framed) BPS states,
arXiv:1301.3065.

\vspace{2mm}

[Dav1] B. Davison, Invariance of orientation data for ind-constructible Calabi-Yau  $A_{\infty}$-categories under derived equivalence, arXiv:1006.5475.

\vspace{2mm}

[Dav2] B. Davison,  The critical CoHA of a self dual quiver with potential, arXiv:1311.7172.
\vspace{2mm}

[Dav3] B. Davison,  Purity of critical cohomology and Kac's conjecture, arXiv:1311.6989.

\vspace{2mm}

[DiHuSo] D.-E. Diaconescu, Z. Hua, Y. Soibelman, HOMFLY polynomials, stable pairs and Donaldson-Thomas invariants,arXiv:1202.4651.

\vspace{2mm}

[DiSVa] D.-E. Diaconescu, V. Shende, C. Vafa, Large N duality, lagrangian cycles, and algebraic knots,
arXiv:1111.6533.
    
\vspace{2mm}

[DyKap] T. Dyckerhoff, M. Kapranov, Higher Segal spaces I, arXiv:1212.3563.

\vspace{2mm}

[Ef] A. Efimov, Cohomological Hall algebra of a symmetric quiver,
 arXiv:1103.2736.
\vspace{2mm}

[EnRe] J. Engel, M. Reineke, Smooth models of quiver moduli, arXiv:0706.43.06.

\vspace{2mm}

[Fra] H. Franzen, On Cohomology Rings of Non-Commutative Hilbert Schemes and CoHa-Module, arXiv:1312.1499.

\vspace{2mm}
[Ga1] D. Gaiotto, ${\mathcal N}=2$ dualities, arXiv:0904.2715.
\vspace{2mm}

[GaMoNe-2] D. Gaiotto, G. Moore, D. Neitzke, Wall-crossing, Hitchin systems and WKB approximations, arXiv:0907.3987.

\vspace{2mm}

[GorORS] E.Gorsky, A. Oblomkov, J. Rasmussen, V. Shende, Torus knots and the rational DAHA,  arXiv:1207.4523.

\vspace{2mm}

[GuSto] S. Gukov, M. Stosic, Homological algebra of knots and BPS states, arXiv:1112.0030. 
   
\vspace{2mm}  

[GuSchVa]  S. Gukov, A. Schwarz, C. Vafa, Khovanov-Rozansky Homology and Topological Strings,
arXiv:hep-th/0412243.
 \vspace{2mm}
[H] D. Huybrechts, Introduction to stability conditions,  arXiv:1111.1745.

\vspace{2mm}

[HaMo1] J. Harvey, G. Moore, Algebras, BPS States, and Strings,  arXiv:hep-th/9510182.

\vspace{2mm}

[HaMo2] J. Harvey, G. Moore, On the algebras of BPS states, arXiv:hep-th/9609017.

\vspace{2mm}

[Hu] Z. Hua. Spin structure on moduli space of sheaves on Calabi-Yau threefold, arXiv:1212.3790.

\vspace{2mm}

[KaKoPa] L. Katzarkov, M. Kontsevich, T. Pantev, Hodge theoretic aspects of mirror symmetry, arXiv:0806.0107.
\vspace{2mm}

[KatzLiu] S. Katz, C. Liu, Enumerative geometry of stable maps with Lagrangian boundary conditions and multiple covers of the disc, arXiv:math/0103074. 
    
\vspace{2mm}

[Ko1] M. Kontsevich, Holonomic $D$-modules and positive characteristic, arXiv:1010.2908.

\vspace{2mm}

[KoSo1] M. Kontsevich, Y. Soibelman, Stability structures, motivic Donaldson-Thomas invariants and cluster transformations, arXiv:0811.2435.

\vspace{2mm}

[KoSo2] M. Kontsevich, Y. Soibelman, Affine structures and non-archimedean analytic spaces, math.AG/0406564.

\vspace{2mm}

[KoSo3] M. Kontsevich, Y. Soibelman, Notes on A-infinity algebras, A-infinity categories and non-commutative geometry. I, math.RA/0606241.

\vspace{2mm}

[KoSo4] M. Kontsevich, Y. Soibelman,Deformations of algebras over operads and Deligne's conjecture, arXiv:math/0001151.

\vspace{2mm}

[KoSo5] M. Kontsevich, Y. Soibelman, Cohomological Hall algebra, exponential Hodge structures and motivic Donaldson-Thomas invariants, arXiv:1006.2706.

\vspace{2mm}

[KoSo6] M. Kontsevich, Y. Soibelman,
Homological mirror symmetry and torus fibrations, arXiv:math/0011041.

\vspace{2mm}

[KoSo7] M. Kontsevich, Y. Soibelman, Lectures on motivic Donaldson-Thomas invariants and wall-crossing formulas, preprint 2010.

\vspace{2mm}

[KoSo8] M. Kontsevich, Y. Soibelman, Wall-crossing structures in Donaldson-Thomas invariants, integrable systems and Mirror Symmetry,  arXiv:1303.3253.

\vspace{2mm}

[Mau1] D. Maulik, Stable pairs and the HOMFLY polynomial, arXiv:1210.6323.

\vspace{2mm}

[MauOk] D. Maulik, A. Okounkov, Quantum Groups and Quantum Cohomology,  arXiv:1211.1287.

\vspace{2mm}

[Nak1] H. Nakajima, AGT conjecture and perverse sheaves on instanton moduli spaces, talk at the Simons Simposium on BPS states and knot invariants, April 2012.

\vspace{2mm}

[Nak2] H. Nakajima, Lectures on Hilbert Schemes of Points on Surfaces, University Lecture Series, AMS, 1999.
\vspace{2mm}

[NagNak] K. Nagao, H. Nakajima, Counting invariant of perverse coherent sheaves and its wall-crossing, arXiv:0809.2992.

\vspace{2mm}

[Neg] A. Negut, The Shuffle Algebra Revisited, arXiv:1209.3349.

\vspace{2mm}

[ORS] A. Oblomkov, J. Rasmussen, V. Shende, The Hilbert scheme of a plane singularity and the HOMFLY homology of its link, arXiv:1201.2115.

\vspace{2mm}

[PaToVaVe] T. Pantev, B. Toen, M. Vaqui\'e, G. Vezzosi, Shifted Symplectic Structures,  arXiv:1111.3209.

\vspace{2mm}

[PT1] R. Pandharipande, R. P. Thomas, Stable pairs and BPS invariants, arXiv:0711.3899.

\vspace{2mm}

[Re1] M. Reineke, Framed quiver moduli, cohomology, and quantum groups, arXiv:math/0411101.

\vspace{2mm}

[Re2] M. Reineke, The Harder-Narasimhan system in quantum groups and cohomology of quiver moduli, arXiv:math/0204059.

\vspace{2mm}

[Rim] R. Rimanyi, On the Cohomological Hall Algebra of Dynkin quivers,  arXiv:1303.3399.

\vspace{2mm}

[Sch1] O. Schiffmann, Lectures on Hall algebras, arXiv:math/0611617.

\vspace{2mm}

[SchV] O. Schiffmann, E. Vasserot, Cherednik algebras, W-algebras and equivariant cohomology of the moduli space of instantons on ${\mathbb A}^2$, arXiv: 1202.2756.

\vspace{2mm} [Se] P. Seidel, Fukaya categories and Picard-Lefschetz theory, Zurich Lectures in Advanced Mathematics, European Mathematical Society, 2008.

\vspace{2mm}

[So1] Y. Soibelman, talks at the Simons Simposium on BPS states and knot invariants, US Virgin Islands, April 2012 and AIM workshop on Donaldson-Thomas invariants and singularity theory, Budapest, May 2012.

\vspace{2mm}

[Sz1] B. Szendroi, Non-commutative Donaldson-Thomas theory and the conifold, arXiv:0705.3419.

\vspace{2mm}

[Sz2] B. Szendroi, Nekrasov's Partition Function and Refined Donaldson-Thomas Theory: the Rank One Case, arXiv:1210.5181.

\vspace{2mm}

[Tod] Y. Toda, Stability conditions and curve counting invariants on Calabi-Yau $3$-folds, arXiv:1103.4229.

\vspace{3mm}

{\it address: Department of Mathematics, KSU, Manhattan, KS 66506, USA

email: {soibel@math.ksu.edu} }

\end{document}